\documentclass[12pt,a4paper]{article}
\newtheorem{theorem}{\bf Theorem}[section]
\newtheorem{lemma}{\bf Lemma}[section]

\newtheorem{question}{\bf Question}[section]
\newtheorem{conjecture}{\bf Conjecture}[section]

\begin{document}
\author{
{\bf Agata Smoktunowicz}\\
{\small\em Maxwell Institute of Mathematical Sciences}\\
{\small\em School of Mathematics, University of Edinburgh}\\
{\small\em Mayfield Road, Edinburgh EH9 3JZ, Scotland, UK}\\
{\small\em A.Smoktunowicz@ed.ac.uk}\\}
\title{Makar-Limanov's conjecture on free subalgebras
\thanks{
 This work was supported by Grant No. EPSRC
 EP/D071674/1.
 }}
\date{ }
\maketitle
\begin{abstract}

It is proved that over every countable field $K$ there is a nil
algebra $R$ such that the algebra obtained from $R$ by extending
the field $K$ contains noncommutative free subalgebras of
arbitrarily high rank.

It is also shown that over every countable field  $K$ there is an
algebra $R$ without noncommutative  free subalgebras of rank two
such that the algebra obtained from $R$ by extending the field $K$
contains a noncommutative  free subalgebra of rank two.
  This answers a question of Makar-Limanov  \cite{lenny}.
  \end{abstract}

\noindent\textit{Mathematics Subject Classification MSC2000:
16S10, 16N40, 16W50, 16U99} \noindent\textit{Key words: free
subalgebras, extensions of algebras, nil rings}
\section{Introduction}
 In the last forty
 years free subobjects in groups and algebras
have been extensively studied by many authors and enormous
progress has been made \cite{anick, chiba, hedges, lorenz, mal,
mliman, rv, salwa, zelmanov}.
  In the influential paper of Makar-Limanov \cite{ml} several interesting
 open questions have been asked. In particular Makar-Limanov  conjectured
 that {\em if $R$ is a finitely generated infinite dimensional
   algebraic division algebra then $R$ contains a free subalgebra in two
   generators.}
  Another question along this line was asked by Anick \cite{anick} in mid 1980's:
    {\em Let R be a finitely presented algebra with exponential
growth. Does it follow that R contains a free subalgebra in
    two generators?}
 In the same paper he shown that finitely presented monomial algebras
    with exponential growth  contain free subalgebras
    in two generators \cite{anick}. In \cite{mal}
Makar-Limanov proved that the quotient algebra of the Weyl algebra
contains a free subalgebra in two generators.
 He also conjectured that the following
 holds.
 \begin{conjecture}[Makar-Limanov, \cite{lenny}, \cite{jb}] If  R is an algebra without free subalgebras of rank
  two and $S$ is an extensions of R obtained
 by extending  the field K then $S$ doesn't contain a free $K$- algebra of rank
 two. \end{conjecture}
Makar-Limanov mentioned that the truth of this conjecture would
imply that we need only to consider algebras over uncountable
fields in his mentioned above conjecture on the division algebras
\cite{ml}.
 Conjecture $1.1$ in the case of skew-fields, as stated in \cite{ml}, attracted a
 lot of attention and is known to be true in several important
 cases \cite{ber, chiba, fgm, fgs, li, mliman, rv, shi}.
 In 1996 Reichstein showed that Conjecture $1.1$ holds for algebras
 over uncountable fields \cite{r}. The purpose of this paper is to show
 that the situation is completely different for algebras over countable
 fields, as shown in the next theorem.
\begin{theorem}
 Over every countable field $K$ there is an algebra $A$ without
 free noncommutative subalgebras of rank two such that the
 polynomial ring $A[x]$ in one indeterminate $x$ over $A$ contains
 a free noncommutative $K$-algebra of rank two.
\end{theorem}
 Note that if an algebra contains a noncommutative free algebra of
 rank two then it also contains a noncommutative free algebra of
 arbitrarily high rank.
 As an application the following result is obtained.
\begin{theorem} For every
countable field $K$ there is a field $F$ with $K\subseteq F$ and a
$K$-algebra $A$ without noncommutative free subalgebras of rank
two such that the algebra $A\otimes_{K}F$ contains a
noncommutative free $K$-subalgebra of rank two.
\end{theorem}
 In the case of skew-fields Makar-Limanov conjecture is still open.

 A ring $R$ is nil if every element $r\in R$ is nilpotent, i.e.
 for every $r\in R$ there is $n$ such that $r^{n}=0$. Jacobson radical
  rings and nil rings are
 useful for investigating the general structure of rings. In
 addition nil rings have applications in group theory. For example
  the famous construction of Golod and Shafarevich, \cite{faith, lam}, in the 1960s produced a finitely
generated nil algebra that was not nilpotent. This was then used
to construct a counterexample to the Burnside Conjecture, one of
the biggest outstanding problems in group theory at that time.
 The Golod-Shafarevich construction gave also a counterexample to the
 Kurosh Problem: {\em let $R$ be a finitely generated
algebra over a field $F$ such that $R$ is algebraic over $F$, is
$R$ finite dimensional over $F$?} However, the Kurosh Problem is
still open for the key special case of a division ring.
 There are connections with problems in nil rings. A nil element
is obviously algebraic, and in the converse direction, it is
possible to construct an associated graded algebra connected with
an algebraic algebra in such a way that the positive part is a
graded nil algebra  \cite{smok1}.

 It was shown by Amitsur in 1973 that if $R$ is a nil algebra over
 an uncountable field then polynomial rings in many commuting variables over $R$ are
also nil \cite{faith, lam}.
 However in general polynomial rings over nil rings need not be nil \cite{smok, smok1}.
 Our next result shows that polynomial rings over some nil rings
 contain noncommutative free algebras of rank two, and hence are very far from being nil.
\begin{theorem}
 Over every countable field $K$ there is a nil algebra $N$ such
 that the polynomial ring $N[X_{1}, \ldots ,X_{6}]$
 in six commuting indeterminates $X_{1}, \ldots ,X_{6}$ over $N$
 contains a noncommutative free $K$-algebra of rank two.
\end{theorem}
 As an application the following result is induced.
\begin{theorem}
 Over every countable field $K$ there is a field $F$, $K\subseteq
 F$ and a nil algebra $R$ such that the algebra $R\otimes_{K}F$ contains a
  noncommutative free $K$-algebra of rank two.
\end{theorem}
\section{Notations}
  Let $K$ be a countable field and let $A$ be the free $K$- algebra
  generated by elements $x_{1}, x_{2}$, $x_{3}$, $y_{1}$, $y_{2}$, $y_{3}$.
    Let $G= \{x_{1}, x_{2}, x_{3}, y_{1}, y_{2}, y_{3}\}$.
   We say that an element $w\in R$ is a monomial, and write $w\in
   M$, if $w$ is a product of elements from $G$. Given $e\in G, w\in M$
   by $\deg_{e}(w)$ we will denote the number of occurrences of
    $e$ in $w$. By $M_{i}$ we denote the set of monomials of degree $i$.
    Let $H_{i}$ be the $K$-linear space spanned by elements from $M_{i}$, i.e.
    $H_{m_{i}}=KM_{i}=span_{K}M_{i}$.
    Let $D$ be the free $K$- algebra generated by
elements $x$, $y$. Denote $x=z_{1}$, $y=z_{2}$. By $P\subseteq D$
we will denote the set of all monomials in $x,y$, and by $P_{i}$
  the set of monomials of degree $i$.
 Let $(i_{1}, \ldots ,i_{m}), (j_{1}, \ldots ,j_{t})$ be integers.
  We say  that $(i_{1}, \ldots ,i_{m})\prec (j_{1}, \ldots ,j_{t})$
   if $(i_{1}, \ldots ,i_{m})$ is smaller than $(j_{1}, \ldots ,j_{t})$
    in the lexicographical ordering, i.e. either
 $i_{1}<j_{1}$ or $i_{1}=j_{1}$
and $i_{2}<j_{2}$, etc.
 Introduce a partial ordering on elements of $P$.
Let $z,z'\in P$ and $z=\prod_{k=1}^{m}z_{i_{k}}$
$z'=\prod_{i=1}^{m'}z_{j_{k}}$ where $i_{k}, j_{k}\in \{1, 2\}$
(recall that $z_{1}=x, z_{2}=y$). We will say that $z\prec z'$ if
$m=m'$ and $(i_{1},\ldots ,i_{m})\prec(j_{1},\ldots ,j_{m})$. Let
$\beta:M\rightarrow P$ be a semigroup homomorphism such that
$\beta (x_{1})=\beta (x_{2})=\beta (x_{3})=x$ and $\beta
(y_{1})=\beta (y_{2})=\beta (y_{3})=y$. Given $z\in P$, define
$S(z)=span_{K}\{w\in M_{\deg {z}}: \beta (w)\prec z\}$,
 $Q(z)=span_{K}\{w\in M_{\deg {z}}: \beta (w)=z\}$.
  Similarly, given $z\in M$, define
 $S(z)=span_{K}\{w\in M_{\deg z}: \beta (w)\prec \beta (z)\}$,
 $Q(z)=span_{K}\{w\in M_{\deg z}: \beta (w)=\beta (z)\}$.
 Given integers $n_{1}, \ldots ,n_{6}$ and a
 monomial
  $w\in P\cup M$, let $w(n_{1}, \ldots , n_{6})=\sum\{v\in Q(w): \deg _{x_{1}}v=n_{1}, \deg_{x_{2}}v=n_{2},
  \deg _{x_{3}}v=n_{3},\deg _{y_{1}}v=n_{4},
  \deg _{y_{2}}v=n_{5}, \deg _{y_{3}}v=n_{6}\}$.
  We put $w(n_{1}, \ldots ,n_{6})=0$
  if either $\deg_{x}w\neq n_{1}+n_{2}+n_{3}$ or $\deg_{y}w\neq
n_{4}+n_{5}+n_{6}$, because
   in this case the sum goes over the empty set.
\begin{lemma} For each $z\in P$ the set $U_{z}=\{z(n_{1}, \ldots ,
n_{6}): 0\leq n_{1}, \ldots, n_{n},$ $\deg_{x}z=n_{1}+n_{2}+n_{3},
\deg_{y}z=n_{4}+n_{5}+n_{6}\}$ is a free basis of a right module
$U_{z}A$. Let $z_{1}, \ldots ,z_{n}\in P_{i}$, for some $i$ and
assume that elements $z_{1}, \ldots ,z_{n}$ are pairwise distinct.
Then the set $T=T_{z_{1}}\bigcup T_{z_{2}}\bigcup \ldots \bigcup
T_{z_{n}}$ is a free basis of a right module $TA$.
\end{lemma}
{\bf Proof.} The proof follows from the fact that $A$ is a free
algebra and elements from $U_{z}$ are linear combinations of
pairwise distinct monomials of the same degree.
\begin{lemma}
 Let $0<p,r$ be natural numbers and let $z=uv$ where $z\in P_{p+r}$, $u\in P_{r}$,
$v\in P_{p}$.
 Then, for arbitrary integers $n_{1}, \ldots , n_{t}$, and
 $r<p+r$ we have
 $z(n_{1}, \ldots, n_{6})=\sum\{u(r_{1}, \ldots ,
 r_{6})v(n_{1}-r_{1}, \ldots ,
 n_{6}-r_{6}):r_{1}+\ldots +r_{6}=r\}$.
\end{lemma}
{\bf Proof.} Observe first that if $p=1$ then
 $z(n_{1}, \ldots, n_{6})=\sum_{i=1}^{6}u_{i}v_{i}$ where $u_{1}=u(n_{1}-1, n_{2}, n_{3}, n_{4}, n_{5},
 n_{6})$, $u_{2}=u(n_{1}, n_{2}-1, n_{3}, n_{4}, n_{5}, n_{6})$, $u_{3}=
u(n_{1}, n_{2}, n_{3}-1,n_{4},n_{5},n_{6})$, $u_{4}= u(n_{1},
n_{2}, n_{3}, n_{4}-1, n_{5}, n_{6})$, $u_{5}= u(n_{1}, n_{2},
n_{3}, n_{4}, n_{5}-1, n_{6})$, $u_{6}= u(n_{1}, n_{2}, n_{3},
n_{4}, n_{5}, n_{6}-1)$ and
 $v_{1}=v(1,0,0,0,0,0)$, $v_{2}=v(0,1,0,0,0,0)$,
$v_{3}=v(0,0,1,0,0,0)$, \ldots , $v_{6}=v(0,0,0,0,0,1)$.
 Note that if $v=x$ then $v_{4}=v_{5}=v_{6}=0$.
  We will prove Lemma $2.2$ by induction on $n$.
  For $n=2$ the result holds because then $r=p=1$. Suppose the result is true
  for some $n>2$. We will show it is true for $n+1$.
  If $n=r+1$ and $p=1$ then the
  result is true by the above observations. If $p>1$ write
  $v=ww'$ for some $w\in P_{p-1}$, $w'\in P_{1}$.

 Then by the case $p=1$ we have
$z(n_{1}, \ldots, n_{6})=\sum_{i=1}^{6}(uw)_{i}w'_{i}$, where
similarly as in the beginning of the proof $(uw)_{1}=uw(n_{1}-1,
n_{2}, n_{3}, n_{4}, n_{5},
 n_{6})$ and $w'_{1}=w'(1,0,0,0,0,0)$, $(uw)_{2}=uw(n_{1},
n_{2}-1, n_{3}, n_{4}, n_{5},
 n_{6})$ and $w'_{1}=w'(0,1,0,0,0,0)$,etc.

By the inductive assumption, $uw(q_{1}, \ldots, q_{6})=
\sum\{u(r_{1}, \ldots ,
 r_{6})w(q_{1}-r_{1}, \ldots ,
 q_{6}-r_{6}):r_{1}+\ldots +r_{6}=r\}$. Now $(uw)_{1}=\sum\{u(r_{1}, \ldots ,
 r_{6})w(n_{1}-1-r_{1}, n_{2}-r_{2}, \ldots ,
 q_{6}-r_{6}):r_{1}+\ldots +r_{6}=r\}$.

Now $uw_{1}w'_{1}=\sum\{u(r_{1}, \ldots ,
 r_{6})w(n_{1}-1-r_{1}, n_{2}-r_{2}, \ldots ,
 q_{6}-r_{6})w'_{1}:r_{1}+\ldots +r_{6}=r\}$.
  Similarly,
$uw_{2}w'_{2}=\sum\{u(r_{1}, \ldots ,
 r_{6})w(n_{1}-r_{1}, n_{2}-r_{2}-1, n_{3}-r_{3}, \ldots ,
 n_{6}-r_{6})w'_{2}:r_{1}+\ldots +r_{6}=r\}$, etc.
 Therefore, $z(n_{1},\ldots ,n_{6})=\sum\{u(r_{1}, \ldots ,
 r_{6})[w(n_{1}-r_{1}-1, n_{2}-r_{2}, \ldots ,
 n_{6}-r_{6})w'_{1}+w(n_{1}-r_{1}, n_{2}-r_{2}-1, \ldots ,
 n_{6}-r_{6})w'_{2}+\ldots +w(n_{1}-r_{1}, n_{2}-r_{2}, \ldots ,
 n_{6}-r_{6}-1)w'_{6}]:r_{1}+\ldots +r_{6}=r\}$.
 Observe that $w(n_{1}-r_{1}-1, n_{2}-r_{2}, \ldots ,
 n_{6}-r_{6})w'_{1}+w(n_{1}-r_{1}, n_{2}-r_{2}-1, \ldots ,
 n_{6}-r_{6})w'_{2}+\ldots +w(n_{1}-r_{1}, n_{2}-r_{2}, \ldots ,
 n_{6}-r_{6}-1)w'_{6}]=ww'(n_{1}-r_{1}, \ldots ,n_{6}-r_{6})$, as
 in the beginning of the proof.
 Therefore,
 $z(n_{1}, \ldots, n_{6})=\sum\{u(r_{1}, \ldots ,
 r_{6})v(n_{1}-r_{1}, \ldots ,
 n_{6}-r_{6}):r_{1}+\ldots +r_{6}=r\}$, as desired.

\begin{lemma} Let $p,q$ be natural numbers.
  Let  $f:H_{p}\rightarrow H_{p}$, $g:H_{q}\rightarrow H_{q}$,
  and $h:H_{p+q}\rightarrow H_{p+q}$ be $K$-linear
  mappings such that for all $w\in M_{p}$, $w'\in M_{q}$, $h(ww')=f(w)g(w')$.
   Let $z\in P_{p+q}$, $z=uv$,
  $u\in P_{p}$, $v\in P_{q}$. If $h(z(n_{1}, \ldots ,n_{6}))\in h(S(z))$ for all  $n_{1}+\ldots
+n_{6}=p+q$ then either $f(u(p_{1}, \ldots, p_{6}))\in f(S(u))$
for all $p_{1}+\ldots +p_{6}=p$ or
     $g(v(q_{1},\ldots , q_{6}))\in g(S(v))$
      for all  $q_{1}+\ldots +q_{6}=q$.
 \end{lemma}
 {\bf Proof.} Suppose that the result does not hold.
  Let $(p_{1}, \ldots, p_{6})$ and $(q_{1}, \ldots , q_{6})$
   be minimal with respect to the ordering $\prec $ and such that
  $p_{1}+ \ldots + p_{6}=p$, $q_{1}+\ldots + q_{6}=q$ and
 $f(u(p_{1}, \ldots, p_{6}))\notin f(S(u))$,
 $g(v(q_{1},\ldots , q_{6}))\notin g(S(v))$.
 Let $D=H_{p}\cap f(S(u))$ and $B=H_{q}\cap g(S(v))$. By Lemma $2.2$,
 $z(p_{1}+q_{1}, \ldots , p_{6}+ q_{6})=\sum_{r_{1}+\ldots+r_{6}=p}u(r_{1},\ldots , r_{6})
 v(p_{1}+q_{1}-r_{1}, \ldots , p_{6}+q_{6}-r_{6})$.
  It follows that $h(z(p_{1}+q_{1}, \ldots , p_{6}+q_{6}))
  =\sum_{r_{1}+\ldots +r_{6}=p}f(u(r_{1},\ldots , r_{6}))
  g(v(p_{1}+q_{1}-r_{1}, \ldots , p_{6}+q_{6}-r_{6}))$.
 Note that if $(p_{1}, \ldots , p_{6})\prec (r_{1}, \ldots , r_{6}) $
  with respect to the lexicographical
 ordering then $(p_{1}+q_{1}-r_{1}, \ldots , p_{6}+q_{6}-r_{6})\prec
 (q_{1}, \ldots , q_{6})$.
  By the assumptions about the minimality of
  $(p_{1}, \ldots , p_{6})$ if $(r_{1}, \ldots , r_{6})\prec(p_{1},\ldots , p_{6})$
   then $f(u(r_{1}, \ldots , r_{6}))\in f(S(u))$.
  Similarly, if $(v_{1}, \ldots , v_{6})\prec (q_{1},
  \ldots , q_{6})$ then $g(v(v_{1}, \ldots , v_{6}))\in g(S((v))$.
 Therefore $h(z(p_{1}+q_{1}, \ldots , p_{6}+q_{6}))\in
 h(z(p_{1},\ldots , p_{6}))g(z(q_{1}, \ldots , q_{6}))+DH_{q}+H_{p}B$.
  By the assumptions of our theorem, $h(z(p_{1}+q_{1}, \ldots , p_{6}+q_{6}))\in
  h(S(z))$. Note that
   since $A$ is generated in degree one $S(z) \subseteq H_{p}S(v)+S(u)H_{q}$ and so
$h(S(z)) \subseteq H_{p}g(S(v))+f(S(u))H_{q}=H_{p}D+BH_{q}$. It
follows that
    $h(z(p_{1}+q_{1}, \ldots , p_{6}+q_{6}))\in DH_{q}+H_{p}B$.
   Therefore, $f(z(p_{1}, \ldots , p_{6}))g(z(q_{1}, \ldots , q_{6}))\in
   DH_{q}+H_{p}B$. Recall that $f(z(p_{1}, \ldots , p_{6}))\in H_{p}$
   and $D\in H_{p}$.
  Therefore either $f(u(p_{1}, \ldots , p_{6}))\in D\subseteq f(S(u))$ or
   $g(v(q_{1}, \ldots , q_{6}))\in B\subseteq g(S(v))$ a contradiction.
\begin{lemma}
 Let $p, r$ be integers such that $p>10^{8}$, $r>10p$, $40$ divides $p+r$. Let
 $f:H_{p}\rightarrow H_{p}$, $g:H_{r+p}\rightarrow H_{r+p}$
be $K$-linear mappings such that for $w\in M_{r}$,
 $w'\in M_{p}$, $g(ww')=wf(w')$. Let $z=uv$, $z\in M_{p+r},u\in M_{r}$, $v\in
 M_{p}$.
 Suppose that for all $n_{1}+\ldots +n_{6}=p+r$, we have
  $$g(z(n_{1}, \ldots , n_{6}))\in \sum_{r_{1}, \ldots
  ,r_{6}:r_{1}+\ldots r_{6}=r}u(r_{1}, \ldots , r_{6})f(S(v))+c
   +\sum_{i=1}^{10^{-4}(r+p)^{2}}Kh_{i}$$
   for some $h_{i}\in H_{p+r}$, and some $c\in \sum_{w}wA$ where $w\in M_{r}$ are monomials  which
    are linearly independent from the elements
   $z(r_{1}, \ldots, r_{6})$ with $r_{1}+\ldots +r_{6}=r$.
    Then $f(v(p_{1}, \ldots, p_{6}))\in f(S(v))$ for all $p_{1}+\ldots +p_{6}=p$.
\end{lemma}
{\bf Proof.} We may assume that $\deg_{x}z\geq {\deg z\over
2}={p+r\over 2}$. In the case when $\deg_{y}z\geq {\deg z\over 2}$
the proof is similar. Note that $f(z(p_{1}, \ldots ,p_{6}))=0$ if
$p_{i}<0$ for some $i$, because then $z(p_{1}, \ldots ,p_{6})=0$.
  Hence, it suffices to show that each $f(v(p_{1}, \ldots
  ,p_{6}))$ is a linear combination of
 $f(v(q_{1}, \ldots , q_{6}))$ with $(q_{1}, \ldots , q_{6})\prec (p_{1}, \ldots , p_{6})$ and elements from $f(S(v))$.
  Let $q_{1}, \ldots ,q_{6}$ be such that $v(q_{1},\ldots
  ,q_{6})\neq 0$. Then
  $\deg_{x}v=q_{1}+q_{2}+q_{3}$ and $\deg_{y}v=q_{4}+q_{5}+q_{6}$ by the
  definition of $v(q_{1}, \ldots ,q_{6})$.
 We will show that $f(v(q_{1}, \ldots ,q_{6}))=0$.
 Let $S=\{(n_{1}\ldots  , n_{6}): {1\over 6}(p+r)<n_{1}<(p+r)({1\over 6}+{1\over 40}),
 {1\over 6}(p+r)<n_{2}<(p+r)({1\over 6}+{1\over 40}), n_{1}+n_{2}+n_{3}=\deg_{x}z$ and moreover
 $n_{4}=q_{4}+\deg_{y}u$, $n_{5}=q_{5}$, $n_{6}=q_{6}\}$.

   First we shall prove that $card(S)\geq (p+r)^{2}10^{-4}$.
   Observe that there are at least $(p+r)40^{-1}-2$ natural
   numbers laying between $(p+r){1\over 6}$ and $(p+r)({1\over 6}+{1\over
   40})$. We can choose $((p+r)(40)^{-1}-2)^{2}$ distinct pairs
   $(n_{1}, n_{2})$ such that ${1\over 6}(p+r)<n_{1}<(p+r)({1\over 6}+{1\over
   40})$ and
 ${1\over 6}(p+r)<n_{2}<(p+r)({1\over 6}+{1\over 40})$. For each
 such pair we can choose a natural number $n_{3}$ such that
 $n_{1}+n_{2}+n_{3}=\deg_{x}z$ and $({1\over 6}-{1\over 20})(p+r)\leq n_{3}$
  because $\deg_{x}z\geq {p+r\over 2}$. Since
 $p+r>10^{8}$, we get that $card(S)\geq ((p+r)(40)^{-1}-2)^{2}>10^{-4}(p+r)^{2}$.

  Hence the assumption of the theorem implies that
  $$\sum_{(n_{1}, \ldots , n_{6})\in S}l_{n_{1}, \ldots ,
  n_{6}}g(z(n_{1}, \ldots, n_{6}))\in\sum_{r_{1}, \ldots
  ,r_{6}:r_{1}+\ldots r_{6}=r}u(r_{1}, \ldots , r_{6})f(S(v))+c,$$ for some $l_{n_{1}, \ldots ,
  n_{6}}\in K$, not all of which are zeros ($c$ is as in the thesis). Let $(j_{1}, \ldots ,
  j_{6})$ be the maximal element in $S$, with respect to $\prec $,
  such that $l_{j_{1}, \ldots , j_{6}}\neq 0$.
   Then $g(z(j_{1}, \ldots , j_{6}))=\sum_{ }k_{n_{1}, \ldots ,
   n_{6}}g(z(n_{1}, \ldots , n_{6})) +q$ where
 the sum runs over all $(n_{1}, \ldots ,
   n_{6})\in S$ with  $z(n_{1},\ldots , n_{6})\prec (j_{1}, \ldots , j_{6})$. Moreover, $q\in \sum_{r_{1}, \ldots
  ,r_{6}:r_{1}+\ldots r_{6}=r}u(r_{1}, \ldots , r_{6})f(S(v))+c$
    for some  $k_{r_{1}, \ldots , r_{6}}\in K$.

 Now $g(v(n_{1}, \ldots , n_{6}))=\sum_{r_{1}+\ldots +r_{6}=r}u(r_{1},
 \ldots , r_{6})f(v(n_{1}-r_{1}, \ldots , n_{6}-r_{6}))$, by Lemma $2.2$.
   Similarly,
  $g(z(j_{1}, \ldots , j_{6}))=\sum_{r_{1}+\ldots
  +r_{6}=r}u(r_{1},\ldots , r_{6})f(v(j_{1}-r_{1},\ldots ,
  j_{6}-r_{6}))$.

   Now substitute these expressions in the equation $$g(z(j_{1}, \ldots , j_{6}))=\sum_{ }k_{n_{1}, \ldots ,
   n_{6}}g(z(n_{1}, \ldots , n_{6})) +q.$$
  We get
$\sum_{r_{1}+\ldots
  +r_{6}=r}u(r_{1},\ldots , r_{6})[f(v(j_{1}-r_{1},\ldots ,
  j_{6}-r_{6}))-\sum_{n_{1}, \ldots ,n_{6}\in S}f(v(n_{1}-r_{1}, \ldots ,
  n_{6}-r_{6})]\in
\sum_{r_{1}+\ldots +r_{6}=r}u(r_{1},\ldots , r_{6})S(v)
  +c$ where
 the sum runs over all $(n_{1}, \ldots ,
   n_{6})\in S$ with  $z(n_{1},\ldots , n_{6})\prec (j_{1}, \ldots , j_{6})$.

  Now, compare the elements starting with
  nonzero  $u(r_{1}, \ldots , r_{6})$ (they are linearly independent by Lemma $2.1$).
   We get the following equations

 $f(z(j_{1}-r_{1}, \ldots , j_{6}-r_{6}))\in \sum k_{n_{1},
 \ldots, n_{6}}f(z(n_{1}-r_{1}, \ldots ,
n_{6}-r_{6}))+ f(S(v))$ where the sum runs over all $(n_{1},
\ldots ,n_{6})\in S$ with  $(n_{1}, \ldots , n_{6})\prec
(j_{1},\ldots , j_{6})$ (provided that $u(r_{1}, \ldots
,r_{6})\neq 0$).
 Consider now elements  $r_{1}=j_{1}-q_{1}$,
$r_{2}=j_{2}-q_{2}$, $r_{3}=j_{3}-q_{3}$ and $r_{4}=\deg_{y}u,$
$r_{5}=r_{6}=0$. We will show that $u(r_{1}, \ldots, r_{6})\neq
0$.
 Observe first that all $r_{i}\geq 0$. It follows because,
 the definition of $S$  and the assumption $r>10p$ imply
   that  $j_{i}>p$ for $i=1, 2, 3$. By the assumptions $q_{1}+q_{2}+q_{3}=\deg_{x}v\leq \deg v=p$.
    Hence for the integers
   $r_{1}=j_{1}-q_{1}$, $r_{2}=j_{2}-q_{2}$, $r_{3}=j_{3}-q_{3}$
   are positive and $r_{1}+r_{2}+r_{3}=(j_{1}+j_{2}+j_{3})-(q_{1}+q_{2}+q_{3})=\deg_{x}z-\deg_{x}v=\deg_{x}u$.
 Observe also that $r_{4}+r_{5}+r_{6}=\deg_{y}u$ as
 required. Hence, $u(r_{1}, \ldots ,u_{6})\neq 0$.
 Therefore, $f(z(q_{1}, \ldots , q_{6}))=f(z(j_{1}-r_{1}, \ldots ,
 j_{6}-r_{6}))\in \sum_{n_{1}, \ldots , n_{6}\prec (j_{1}, \ldots ,
 j_{6})}k_{n_{1}, \ldots , n_{6}}f(z(n_{1}-r_{1},\ldots ,
 n_{6}-r_{6}))+f(S(v)).$
  Clearly,
  $(n_{1}-r_{1},\ldots , n_{6}-r_{6})\prec (j_{1}-r_{1},
  j_{2}-r_{2}, j_{6}-r_{6})$, so the result holds.
\section{Some results from other papers}
 In this section we quote some results from \cite{smok}. These results will be used
in the last section to get the main result.

Let $A$ be a $K$- algebra generated by elements $x_{1}, x_{2}$,
$x_{3}$, $y_{1}$, $y_{2}$, $y_{3}$ with gradation one. Write
$A=H_{1}+H_{2}+\ldots $. Recall that $H_{i}=KM_{i}$. We will write
$M_{0}=\{1\}\subseteq K$, $H_{0}=K$.
 Given a number $n$ and a set $F\subseteq A$ by $B_{n}(F)$ we will
 denote the right ideal in $A$ generated by the set
 $\bigcup_{k=0}^{\infty }M_{nk}F$, i.e.,
 $B_{n}(F)=\sum_{k=0}^{\infty }H_{nk}FA$.
 \begin{theorem}
 Let $f_{i}$, $i=1,2, \ldots $ be polynomials in $A$ with
  degrees $t_{i}$, and let $m_{i}$, $i=1,2,
 \ldots $ be an increasing sequence of natural numbers such that
 $m_{i}>6^{6t_{i}}$ and $m_{1}>10^{8}$.  There exists subsets $F_{i}\subseteq
 H_{m_{i}}$ with
 $card (F_{i})<10^{-4}m_{i}^{2}$ such that the ideal $I$ of
 $A$ generated by $f_{i}^{10m_{i+1}}$, $i=1, 2, \ldots $ is
 contained in the right ideal $\sum_{i=0}^{\infty }B_{m_{i+1}}(F_{i})$.
 Moreover, for every $k$, $I\cap H_{m_{k+1}}\subseteq \sum_{i=0}^{k}B_{m_{i+1}}(F_{i})$.
\end{theorem}
{\bf Proof.} Let $I_{i}$ be the smallest homogeneous ideal in $A$
 containing $f_{i}^{10m_{i+1}}$, for $i=1, 2, \ldots $.
By considering algebras generated by $6$ elements instead of $3$
elements and using  the same proof as the proof of Theorem $2$ in
\cite{smok} for $k=m_{i}$, $w=m_{i+1}$, $f=f_{i}$ and changing
constants from $3$ to $6$, we get the following result. There
exists a set $F_{i}\subseteq H_{m_{i}}$, such that
card$F_{i}<m_{i}6^{6t_{i}}t_{i}^{2}$ such that the (two sided)
ideal of $A$ generated by $f_{i}^{10m_{i+1}}$ is contained in
$B_{m_{i+1}}(F_{i})$.
 Note that card$F_{i}<10^{-4}m_{i}^{2}$ since
 $m_{i}>6^{6t_{i}}$ and $m_{i}>m_{1}>10^{8}$ by the assumptions.
 Observe now that $I\subseteq \sum_{i=1}^{\infty }I_{i}$.
  Note that $I_{k+1}$ is generated by elements with degrees larger
  than $m_{k+1}$. Recall that ideals $I_{i}$ are homogeneous.
  Therefore, $I\cap H_{m_{k+1}}\subseteq \sum_{i=1}^{k}I_{i}$.
 Hence, $I\cap H_{m_{k+1}}\subseteq
 \sum_{i=1}^{k}B_{m_{i+1}}(F_{i})$ as required.
 This finishes the proof.

 Let mappings $R_{i}:H_{m_i}\rightarrow H_{m_i}$ and $c_{R_{i}(F_{i})}$  be defined as in
 section $2$ in \cite{smok}
 with $F_{i}=\{f_{i,1}, \ldots ,f_{i,r_{i}}\}\subseteq H_{m_{i}}$ be as in Theorem $3.1$.
 Recall that $c_{R_{i}(F_{i})}:H_{m_{i}}\rightarrow H_{m_{i}}$ is
 a $K$-linear mapping with ker$c_{R_{i}(F_{i})}=\{R_{i}(f_{i,1}),
 \ldots ,R_{i}(f_{i,r_{i}})\}$. Given $w=x_{1}\ldots
 x_{m_{i+1}}\in M_{m_{i+1}}$, $R_{i+1}:H_{m_{i+1}}\rightarrow H_{m_{i+1}}$ is a
 $K$-linear mapping such that
  $$R_{i+1}(w)=c_{R_{i}(F_{i})}
  (R_{i}(x_{1}\ldots
  x_{m_{i}}))\prod_{j=2}^{m_{i+1}m_{i}^{-1}}
  R_{i}(x_{(j-1)m_{i}+1}\ldots
  x_{jm_{i}}).$$

 Moreover, $R_{1}=Id$.
 The fact that the algebra $A$ is generated by $6$ elements
 instead of $3$ elements doesn't change the proof of Theorem $4$
 in \cite{smok}.
\begin{theorem}[Theorem $4$, \cite{smok}]
 Suppose that $w\in H_{m_{l+1}}\cap \sum_{i=0}^{l}B_{m_{i+1}}(F_{i})$. Then
 $R_{l+1}(w)=0$.
\end{theorem}
\section{Linear mappings}
 In this section we will prove some technical results about  the
 mappings $R_{i}$. The algebra $A=H_{1}+H_{2}+\ldots $ is as in
 the previous sections.  We will
 use the following notations. $M_{0}=\{1\}$ and $H_{0}=K$.
 In this section we will assume that $R_{i}:H_{m_{i}}\rightarrow
H_{m_{i}}$ are as in section $3$ and moreover  $40m_{i}$ divides
$m_{i+1}$ and $m_{i+1}>2^{i+101}m_{i}$, $m_{1}>10^{8}$ for $i=1,
2, \ldots $.
\begin{lemma} Let $k$ be a natural number. Then
 there are non-negative integers $e_{i}, d_{i}$ with
$\sum_{i}e_{i}>50\sum_{i}d_{i}$ and $\sum_{i}e_{i}+d_{i}=m_{k}$
 such that if $w\in M_{m_{i}}$ and $w=\prod_{i}u_{i}v_{i}$ with
 $u_{i}\in M_{e_{i}}$, $v_{i}\in M_{d_{i}}$ then
 $R_{k}(w)=\prod_{i}u_{i}g_{i,k}(v_{i})$ for some
  $K$-linear mappings $g_{i,k}:H_{d_{i}}\rightarrow  H_{d_{i}}$.

Let $\sigma $ be a permutation on a set of $m_{k}$ elements, such
that $(\prod_{i=1} u_{i}v_{i})^{\sigma
}=\prod_{i}u_{i}\prod_{i}v_{i}$.
 Denote $u=\prod_{i}u_{i}$, $v=\prod_{i}v_{i}$.
  Let $T_{k}(uv)=R_{k}((uv)^{\sigma
^{-1}})^{\sigma }$. Then $T_{k}(uv)=uf{k}(v)$, where
$f_{k}:H_{\deg v}\rightarrow H_{\deg v}$ is a $K$-linear mapping
defined as follows
$f_{k}(v)=f_{k}(\prod_{i}v_{i})=\prod_{i}g_{i,k}(v_{i})$.
\end{lemma}
{\bf Proof.}  The proof of the first part of Lemma $4.1$ is the
same as the proof of Theorem $6$ in \cite{smok}. Note that
$e_{1}=0$ and $u_{1}=1\in K$.
 To prove the second part of Lemma $4.1$, observe that $T_{k}(uv)=R_{k}(w)^{\sigma
}=R_{k}(\prod_{i}u_{i}v_{i})^{\sigma
}=(\prod_{i}u_{i}g_{i,k}(v_{i}))^{\sigma
}=\prod_{i}u_{i}\prod_{i}g_{i,k}(v_{i})= uf_{k}(v)$, as required.
\begin{lemma} Let $w=\prod_{i}u_{i}v_{i}$,
 $u=\prod_{i}u_{i}$, $v=\prod_{i}v_{i}$, $e_{i}$, $d_{i}$, $T_{k}$ be as
in Lemma $4.1$. Let $k$ be a natural number. Then
$$(R_{k}(S(w)))^{\sigma }\subseteq \sum_{c\in M_{\deg u}: c\notin Q(u)}cA
+\sum_{c\in M: c\in Q(u)}cf_{k}(S(v)).$$
  Moreover $$R(w(n_{1}, \ldots , n_{6}))=(\sum_{p_{1}+\ldots
+p_{6}=\deg u}u(p_{1}\ldots p_{6})f_{k}(v(n_{1}-p_{1}, \ldots
,n_{6}-p_{6})))^{\sigma ^{-1}}, $$ for all $n_{1}, \ldots ,n_{6}$.
\end{lemma} {\bf Proof.} Observe first that  $S(w)$ is a linear combination of some
elements $t=\prod_{i}q_{i}r_{i}$ with $q_{i}\in M_{e_{i}}$,
$r_{i}\in M_{d_{i}}$. If $\prod_{i}q_{i}\in Q(u)$ then $q_{i}\in
Q(u_{i})$ for each $i$. In this case, since
$\prod_{i}q_{i}r_{i}\in S(w)$ we have $\prod_{i}r_{i}\in S(v)$.

By the definition of the mapping $R_{k}$ we have
$R_{k}(t)=\prod_{i}q_{i}g_{i,k}(r_{i})$. Now $(R_{k}(t))^{\sigma
}=\prod_{i}q_{i}\prod_{i}g_{i,k}(r_{i})=\prod_{i}q_{i}f_{k}(\prod_{i}r_{i})$.
  Recall that, if $\prod_{i}q_{i}\in Q(u)$ then $\prod_{i}r_{i}\in
  S(v)$.
  Consequently, $f_{k}(\prod_{i}r_{i})\in f_{k}(S(v))$, and so
  $(R_{k}(S(w)))^{\sigma }\subseteq \sum_{c\in M_{\deg u}: c\notin Q(u)}cA+
  \sum_{c\in M: c\in Q(u)}cf_{k}(S(v)).$

We will now prove the second part of the theorem.  Let $z=uv$, by
Lemma $2.2$, we have $\sum_{p_{1}+\ldots +p_{6}=\deg u}
 u(p_{1}\ldots p_{6})f_{k}(v(n_{1}-p_{1}, \ldots
,n_{6}-p_{6}))=T_{k}(z(n_{1}, \ldots ,n_{6}))$. Note that
$z^{\sigma ^{-1}}=w$. Therefore, $T_{k}(z(n_{1}, \ldots ,n_{6}))=
R_{k}(z(n_{1}, \ldots ,n_{6})^{\sigma ^{-1}})^{\sigma}=
R_{k}(w(n_{1}, \ldots ,n_{6}))^{\sigma}.$ The result follows.
 \begin{lemma}  Let $w=\prod_{i}u_{i}v_{i}$,
  $u=\prod_{i}u_{i}$, $v=\prod_{i}v_{i}$, $e_{i}$, $d_{i}$, $T_{k}$, $f_{k}$ be as
in Lemma $4.2$. Let $k$ be a natural number.
  Suppose that $f_{k}(v(n_{1},\ldots ,n_{6})\in f_{k}(S(v))$ for all $n_{1}+\ldots +n_{6}=\deg v$.
 Then $R_{k}(w(n_{1}, \ldots ,n_{6}))\subseteq  R_{k}(S(w))$ for all
$n_{1}+\ldots +n_{6}=m_{i}$.
\end{lemma}
 {\bf Proof.} By the assumption that $f_{k}(v(n_{i}, \ldots
,n_{6})\in f_{k}(S(v))$. Let $z=uv$. Hence, by Lemma $2.2$,
$z(n_{1}, \ldots ,n_{6})\in Q(u)S(v)$ for all $n_{1}, \ldots
,n_{6}$. Consequently, $T_{k}(z(n_{1}, \ldots ,n_{6}))\in
Q(u)f_{k}(S(v))$ for all $n_{1}, \ldots ,n_{6}$.
 Now, by Lemma $4.1$ we have $R_{k}(w(n_{1},
\ldots ,n_{6})) \in [Q(u)S(v)]^{\sigma ^{-1}}$. An element in
$S(v)$ is a linear combination of some elements $\prod_{i}r_{i}\in
S(v)$, with $r_{i}\in M_{d_{i}}$.
 An element $p\in Q(u)$ is a linear combination
 of products $\prod_{i} q_{i}$, with $q_{i}\in Q(u_{i})$.
 Therefore elements from the set
  $Q(u)S(v)$ are linear combinations of products
$\prod_{i}q_{i}\prod_{i}r_{i}$. It follows that  elements from the
set $[Q(u)f_{k}(S(v))]^{\sigma ^{-1}}$ are linear combinations of
products $[\prod_{i}q_{i}\prod_{i}g_{i,k}(r_{i})]^{\sigma ^{-1}}=
\prod_{i}q_{i}g_{i,k}(r_{i})=R_{k}(\prod_{i}q_{i}r_{i})$.
 It follows that $\prod_{i}q_{i}r_{i}\in S(w)$ since
 $\prod_{i}q_{i}\in Q(u)$ and $\prod_{i}r_{i}\in S(v)$, as required.
\begin{theorem} Let $T_{k}$, $u=\prod_{i}u_{i}$, $v=\prod_{i}v_{i}$,
  $w=\prod_{i}u_{i}v_{i}$, be as in Lemma $4.2$. If $R_{k}(w(n_{1}, \ldots ,n_{6}))\subseteq
R_{k}(S(w))+\sum_{i=1}^{m_{k}^{2}10^{-4}}Kg_{i}$ for some
$g_{i}\in A$ then $R_{k}(w(n_{1}, \ldots ,n_{6}))\subseteq
R_{k}(S(w))$ for all $n_{1}, \ldots ,n_{6}$.
\end{theorem}
{\bf Proof.}  By Lemma $4.2$ we have $T_{k}(z(n_{1},\ldots ,
n_{6}))= R_{k}(\bar {z}(n_{1}, \ldots , n_{6})^{\sigma
^{-1}})^{\sigma }=R_{k}(w(n_{1}, \ldots ,n_{6}))^{\sigma }$ for
all $n_{1}, \ldots ,n_{6}$.
 By assumption $R_{k}(w(n_{1}, \ldots
,n_{6}))^{\sigma }\subseteq
R_{k}(S(w))^{\sigma}+\sum_{i=1}^{m_{k}^{2}10^{-4}}Kg_{i}^{\sigma}.$
Denote $g_{i}^{\sigma }=h_{i}$. By Lemma $4.2$
$(R_{k}(S(w)))^{\sigma }\subseteq \sum_{c\in M_{\deg u}: c\notin
Q(u)}cA+\sum_{c\in M:
 c\in Q(u)}cf_{k}(S(v)).$
 It follows that $T_{k}(z(n_{1},\ldots ,n_{6}))\subseteq
 \sum_{c\in M_{\deg u}: c\notin Q(u)}cA+\sum_{c\in M: c\in Q(u)}cf_{k}(S(v))
+\sum_{i=1}^{m_{i}^{2}10^{-4}}Kf_{i}.$
 Therefore $T_{i}$ satisfies
the assumptions of Lemma $2.4$. Consequently, $f_{k}(v(n_{1},
\ldots ,n_{6}))\in f_{k}(S(v))$ for all $n_{1}, \ldots ,n_{6}$.
 By Lemma $4.3$ we get that
$R_{k}(w(n_{1}, \ldots ,n_{6}))\subseteq R_{k}(S(w))$ for all
$n_{1}, \ldots ,n_{6}$, as required.
\begin{theorem}
  Let $i>0$, $F_{i}=\{f_{i,1}, \ldots ,f_{i, r_{i}}\}\subseteq
  H_{m_{i}}$, with $r_{i}<10^{-4}m_{i}^{2}$. For every  monomial $w\in P$ of degree
$m_{i}$ for some $i$,
   there are $n_{1}, \ldots ,n_{6}$ such $n_{1}+\ldots +n_{6}=m_{i}$
   such that  $R_{i}(w(n_{1}, \ldots , n_{6}))\notin R_{i}(S(w))$.
\end{theorem}
{\bf Proof.} Suppose on the contrary.  Let i be the minimal number
such that there is a monomial $w\in P_{m_{i}}$ with
$R_{i}(w(n_{1}, \ldots , n_{6}))\in R_{i}(S(w))$ for all $n_{1},
\ldots ,n_{6}$.
  Clearly $i>1$, since $R_{1}=Id$, $m_{1}>10^{8}$ and $A$ is a free algebra.
 Write $w=w_{1}w_{2}\ldots w_{m_{i}\over m_{i-1}}$ where all
 $w_{i}\in H_{m_{i-1}}$.
  By the definition of $R_{i}$
  and by Lemma $2.3$ we get that either for some $j>1$ we have
  $R_{i-1}(w_{j}(p_{1}, \ldots , p_{6}))\in S(w(j))$  for
   all $p_{1}+ \ldots + p_{6}=m_{i-1}$
 or we have $c_{R_{i-1}(F_{i-1})}(R_{i-1}(w_{1}(p_{1}, \ldots, p_{6}))\in
  c_{R_{i-1}(F_{i-1})}(S(w_{1}))$
for all $p_{1}+\ldots +p_{6}=m_{i-1}$.
 Note that $i$ was minimal,
 and hence the former is impossible. Thus suppose the
 later holds. Then, by the definition of the mapping
  $c_{R_{i-1}(F_{i-1})}$ we have
   $R_{i-1}(w_{1}(n_{1}, \ldots , n_{6})-q(n_{1}, \ldots ,n_{6}))\in
  +\sum_{j=1}^{r_{i-1}}KR_{i-1}(f_{i-1,j})$, for
  some $q(n_{1}, \ldots ,n_{6})\in S(w_{1})$.
  Therefore, $R_{i-1}(w_{1}(n_{1}, \ldots , n_{6}))\in
 R_{i-1}(S(w_{1})) +\sum_{j=1}^{r_{i-1}}KR_{i-1}(f_{i-1,j})$.
 By assumption $r_{i-1}<10^{-4}m_{i-1}^{2}$.
  Theorem $4.1$ applied for $k=i-1$ yields, $R_{i-1}(w_{1}(n_{1}, \ldots , n_{6}))\in
 R_{i-1}(S(w_{1}))$. It is a contradiction, because $i$ was
 minimal.\\
\section{ The main results}
 In this section we will prove Theorems $1.1-1.4$.
 The general idea of the proof of Theorem $1.3$ is a little similar to the proof
 that polynomial rings over nil rings need not be nil, in
 \cite{smok}. Theorems $1.1$, $1.2$ and $1.4$ are consequences of
 Theorem $1.3$.
\\
{\bf Proof of Theorem $1.3$.} Let $K$ be a countable field and let
$A$ be the free noncommutative associative $K$ algebra in
generators $x_{1}, x_{2}, x_{3}, y_{1}, y_{2}, y_{3}$. The field
$K$ is countable so elements of $A$ can be enumerated, say $f_{1},
f_{2}, \ldots $ where degree of $f_{i}$ is $t_{i}$.
 Let $I$ be an ideal in $A$ generated by the homogeneous components
 of elements $f_{i}^{10m_{i+1}}$, $i=1,2, \ldots $ where  $m_{i}$, $i=1,2,
 \ldots $ is an increasing sequence of natural numbers such that
 Let
$40m_{i}$ divide $m_{i+1}$ and $m_{i+1}>2^{i+101}m_{i}$,
$m_{1}>10^{8}$ for $i=1, 2, \ldots $.
  Denote $N=A/I$. Observe that $N$ is nil.
  Let $B$ be the subalgebra of $N[X_{1}, \ldots ,X_{6}]$ generated by elements
   $X=x_{1}X_{1}+x_{2}X_{2}+x_{3}X_{3}+I[X_{1}, \ldots ,X_{6}]$ and element
   $Y=y_{1}X_{4}+y_{2}X_{5}+y_{3}X_{6}+I[X_{1}, \ldots ,X_{6}]$.
 Let $Q$ be the subgroup of $N$ generated by elements $X, Y$ and
 let $P$ be the free subgroup generated by elements $x, y$ as
 in section $2$ and let $\xi:P\rightarrow Q$ be a subgroup
  homomorphism such that $\xi (x)=X$, $\xi (y)=Y$.
We will show that $B$ is a free algebra. Note that the ideal $I$
is homogeneous, hence we only need to show that  linear
combinations of non-zero elements of the same degree are non-zero
(or else all coefficients are zero). Suppose on the contrary. Then
there is $v\in P_{m_{k}} $ for some $k$ such that $\xi (w)\in
\sum_{v\prec w}K\xi (v)$. By rewriting this and comparing elements
 with a pre-fix
 $x_{1}^{n_{1}}x_{2}^{n_{2}}x_{3}^{n_{3}}y_{1}^{n_{4}}y_{2}^{n_{5}}y_{3}^{n_{6}}$
we get that
 $w(n_{1}, \ldots ,n_{6})+I\subseteq S(w)+I$, for all $n_{1}, \ldots , n_{6}$. Therefore,
$w(n_{1}, \ldots ,n_{6})\subseteq S(w)+I$. Note that $w(n_{1},
\ldots ,n_{6})\in H_{\deg w}=H_{m_{k}}$.
  By Theorem $3.1$ there exists subsets $F_{i}\subseteq
 H_{m_{i}}\subseteq A$, with
 $card (F_{i})<10^{-4}m_{i}^{2}$ such that $I\cap H_{m_{k}}\subseteq\sum_{i=1}^{k-1 }B_{m_{i+1}}(F_{i})$.
 It follows that, $w(n_{1}, \ldots ,n_{6})\subseteq S(w)+\sum_{i=1}^{k-1 }B_{m_{i+1}}(F_{i})\cap H_{m_{k}}$.
 By Theorem $3.2$ $R_{k}(\sum_{i=1}^{k-1 }B_{m_{i+1}}(F_{i})\cap
 H_{m_{k}})=0$. Hence,
$R_{k}(w(n_{1}, \ldots ,n_{6}))\subseteq R_{i}(S(w))$, for all
 $n_{1}, \ldots ,n_{6}$. By Theorem $4.2$ it is impossible.
\\
{\bf Proof of Theorem $1.3$.} It follows from Theorem $1.3$ when
we take $F=K\{X_{1}, \ldots ,X_{6}\}$, the field of rational
functions in $6$ commuting indeterminates over $A$ where $A$ is as
in Theorem $1.3$.

{\bf Proof of Theorem $1.1$.} Let $A$ be as in Theorem $1.3$.
Consider rings $R_{0}=A$, $R_{1}=A[X_{1}], R_{2}=A[X_{1}, X_{2}]$,
 $\ldots $ , $R_{6}=A[X_{1}, \ldots ,X_{6}]$.
 Note that $R_{0}$ doesn't contain free algebras of rank two and
 $R_{6}$ contains a free algebra of rank $2$.
  Then there is $0\leq i<6$,
 such that $R_{i}$ doesn't contain
free algebras of rank two and
 $R_{i+1}$ contains a free algebra of rank $6$. Then $R_{i}$
 satisfies the thesis of Theorem $1.1$.

{\bf Proof of Theorem $1.2$.}
 It follows from Theorem $1.1$ when
we take $F=K\{X_{1}, \ldots ,X_{6}\}$, the field of rational
functions in $6$ commuting indeterminates over $A$ where $A$ is as
in Theorem $1.1$.

 {\bf Acknowledgements} The author
would like to thank Jason Bell and Lenny Makar-Limanov for
bringing Conjecture $1.1$ to her attention and to Zinovy
Reichstein for many helpful remarks.


\begin{thebibliography}{27}
  \bibitem{anick} David J Anick, On monomial algebras of
  finite global dimension, Trans. Amer. Math.Soc. {\bf 291} (1985), no. 1, 291-310.
\bibitem{jb} Jason Bell, provate communication, May 2006.
\bibitem{ber} Jason Bell, D Rogalski, Free subalgebras of division
algebras, preprint.
\bibitem{chiba} Katsuo Chiba, Free subgroups and free subsemigroups of
division rings.
 J. Algebra 184 (1996), no. 2, 570--574.
\bibitem{fgm} V Ferreira, J Gonçalves, Z Jairo, A Mandel,
 Arnaldo Free symmetric and unitary pairs in division rings with
 involution.  Internat. J. Algebra Comput.  15  (2005),  no. 1, 15--36.
 \bibitem{fgs} L M Figuerido, J Z Goncalves, M Shirivani,
   Free group algebras in certain division rings, J. Algebra
  {\bf 185} (1996), no. 2, 298--311.
\bibitem{faith} Faith, C., Rings and Things and a Fine Array of Twentieth
Century Associative Algebra. Surveys of The American Math. Soc.,
Vol 65, Providence, 1999, 2nd ed., 2004.
 \bibitem{hedges} M Hedges, The Freiheitssatz for graded
 algebras, J. London Math.Soc. (2) {\bf 35} (1987), no. 3,
 395--405.
 \bibitem{lam} Lam, T., Y., A first course in Noncommutative rings. Second edition,
Graduate texts in Mathematics, 131, Springer-Verlag, New York,
2001.
 \bibitem{li} A I Lichtman, Free subalgebras in division
 rings generated by universal enveloping algebras, Algebra
 Colloq. {\bf 6} (1999), no.2, 15-153.
\bibitem{lorenz} Martin Lorenz, On free subalgebras of
certain division algebras, Proc. Amer. Math. Soc. {\bf 98} (1986),
no. 3, 401--405.
\bibitem{ml} Lenny Makar-Limanov, On free subobjects of skew
fields, Methods in Ring Theory (Antwerp 1983), 281--285. NATO Adv.
Sci. Inst. Ser. C. Math. Phys. Sci. 129, Reidel, Dordrecht, 1984.
\bibitem{mal} Lenny Makar-Limanov, The skew field of
fractions of the Weyl algebra contains a free noncommutative
subalgebra, Comm. Algebra {\bf 11} (1983), no. 17, 2003--2006.
\bibitem{mliman} Lenny Makar-Limanov, On free subsemigroups
of skew fields, Proc. Amer. Math.Soc. {\bf 91} (1984), no. 2,
189-191.
 \bibitem{lenny} Lenny Makar-Limanov, private communication,
 Beijing, June 2007.

\bibitem{r} Zinovy Reichstein, On a qustion of Makar-Limanov,
Proc. Amer.Math.Soc. {\bf 124} (1996), no. 1, 17--19.
\bibitem{rv} Zinovy Reichstein, Nikolaus Vonessen, Free subgroups in division
algebras, Comm.Algebra {\bf 23} (1995), no. 6, 2181--2185.
 \bibitem{salwa} Arkadiusz Salwa, On free subgroups of units
 of rings, Proc. Amer. Math.Soc. {\bf 127} (1999), no. 9,
 2569--2572.
 \bibitem{shi} M Shirvani, J Z Gonçalves,
  Large free algebras in the ring of fractions of skew polynomial
 rings, J. London Math. Soc. (2) 60 (1999), no. 2, 481--489.
\bibitem{smok} Agata Smoktunowicz,
Polynomial rings over nil rings need not be nil, Journal of
Algebra, {\bf 233} (2000), no.2, 427--436.
\bibitem{smok1} Agata Smoktunowicz, Some results in Noncommutative Ring
Theory, Proceedings of the International Congress of
Mathematicians, Madrid, Spain 2006, Vol 2, Invited talks, 2006,
259-269.
\bibitem{zelmanov} Efim
Zelmanov, On groups satysfying
 Golod-Shafarevich condition, New horizons in pro-p groups,
 223--232, Progr. Math., 184, Brikh{\. a}user Boston, Boston, MA
 2000.
\end{thebibliography}
\end{document}